%% file: RobustPhaseBalancing.tex
\documentclass[10pt]{article}
\input{hicss51-packages.tex}
\usepackage[normalem]{ulem}
\usepackage{bbm, xcolor}

\newcommand{\se}[1]{{\textcolor{black}{{#1}}}}

\usepackage{mathtools}

\newcommand{\minimize}{\mathop{\textup{minimize}}}
\newcommand{\maximize}{\mathop{\textup{maximize}}}

\title{Robust Look-ahead Three-phase Balancing of Uncertain Distribution Loads}

\author{Xinbo Geng \\
  Department of Electrical and\\ Computer Engineering\\
  Texas A\&M University \\
  College Station, TX \\
  {\underline{xbgeng@tamu.edu }} \\\And
  Swati Gupta \\
  School of Industrial and\\ Systems Engineering\\ Georgia Institute of Technology\\ Atlanta, GA \\
  {\underline{ swatig@gatech.edu} }\\\And 
  Le Xie \\
  Department of Electrical and\\ Computer Engineering\\
  Texas A\&M University \\
  College Station, TX \\
  {\underline{le.xie@tamu.edu}} \\}

\date{}

\begin{document}
\maketitle
\begin{abstract}
Increasing penetration of highly variable components such as solar generation and electric vehicle charging loads pose significant challenges to keeping three-phase loads balanced in modern distribution systems. Failure to \se{maintain} balance across three phases would lead to asset deterioration and increasing delivery losses. Motivated by the real-world needs to automate and optimize the three-phase balancing decision making, this paper introduces a robust look-ahead optimization framework that \se{pursues} balanced phases in the presence of \se{demand-side} 
uncertainties. \se{We show that} look-ahead moving window optimization \se{can} reduce imbalances among phases at the cost of a \se{limited} number of phase swapping operations. Case studies quantify the improvements of the proposed methods compared with conventional deterministic phase balancing. Discussions on possible benefits of the proposed methods and extensions are presented.
\end{abstract}

\section{Introduction}
Increasing levels of distributed energy resources, together with more active participation of demand side programs, have introduced higher levels of uncertainties to  distribution grid operations. \se{One fundamental task} for distribution system operators \se{(DSOs)} is to keep three phases \se{as balanced as possible} over a long period of time. However, the increasing variability coming from end users \se{requires} \se{DSOs} to revisit this old problem with modern techniques.  \se{Imbalanced three phases} could lead to higher risks of equipment failures \cite{zhu_phase_1998}, increased \se{delivery} losses \cite{ochoa_evaluation_2005}, potential relay malfunctioning \cite{chen_optimal_1999}, additional asset reinforcement costs \cite{ma_quantification_2016}, and issues related with voltage imbalances \cite{yan_investigation_2013,shahnia_voltage_2014,modarresi2017robust,xu_data-driven_2018}. In particular, there is increasing need to develop solutions that can keep \se{three phases} balanced in the presence of high uncertainties from end users over a period of time (e.g. over the course of a day). In this paper, we provide a novel and scalable solution for addressing this problem.

From a \se{DSO}'s perspective, there are three levels of \se{decisions} that \se{can} be made \se{to} \se{ensure reliable and efficient delivery} \se{of electricity to end-users} during normal conditions. At the highest level, it can engage with transmission-level voltage/reactive power optimization routine to regulate its voltage level at the point of interconnection with the backbone grid \cite{geng_voltage_2017,xinbo_geng_chance-constrained_2018}. At the medium level, modern distribution operator could control various sectionalizers and tie switches in order to optimize the \emph{topology} of a distribution system \cite{khodabakhsh_submodular_2017,civanlar_distribution_1988}. At the lowest level, \se{DSOs need} to optimize the assignment of each load (or each cluster of \se{loads}) to appropriate \se{phases} in order to keep the three phase \se{balanced} during a wide range of operating conditions. This paper addresses the issue at the lowest level.

There is a large body of \se{literature} that addresses the issue of keeping \se{three phases balanced} in distribution systems. The phase balancing problem has been traditionally formulated as a mixed integer linear program (MILP) \cite{zhu_phase_1998}. Due to the computational intractability of mixed integer programs, many \se{optimization techniques and heuristics} have been applied \se{to} phase balancing: simulated annealing \cite{zhu_phase_1999}, expert systems \cite{lin_expert_2008}, particle swarm optimization\cite{hooshmand_fuzzy_2012}, immune algorithm \cite{chen_optimal_1999,lin_optimal_2007} and dynamic programming \cite{wang_phase_2013}. These works \cite{zhu_phase_1998,lin_expert_2008,hooshmand_fuzzy_2012,chen_optimal_1999,lin_optimal_2007,wang_phase_2013} typically consider either a single snapshot or use average loads over a long \se{period} of time. \se{In} \cite{dilek_simultaneous_2001,lin_heuristic_2005}, \se{the authors} demonstrate the benefits of extending the phase balancing problem to multiple snapshots and utilizing daily load patterns. It is worth mentioning that \cite{zhu_phase_1998,lin_expert_2008,hooshmand_fuzzy_2012,chen_optimal_1999,lin_optimal_2007,wang_phase_2013,dilek_simultaneous_2001,lin_heuristic_2005} solve the deterministic phase balancing problem. \emph{Uncertainties as well as inter-temporal variabilities} have not been taken into account in the problem of phase balancing. This is the key gap we attempt to bridge in this work.

The reminder of this paper is organized as follows: Section \ref{sec:robust_optimization} introduces robust optimization; Section \ref{sec:phase_balancing} first reviews the deterministic phase balancing problem, which is enhanced to a robust optimization problem in Section \ref{sub:robust_phase_balancing}. The proposed robust look-ahead phase balancing problem is in Section \ref{sub:robust_look_ahead_phase_balancing}. Case studies and discussions are presented in Section \ref{sec:case_study} and \ref{sec:discussions}. Conclusions and future works are in Section \ref{sec:concluding_remarks}.

\section{Robust Optimization: Preliminaries} 
\label{sec:robust_optimization}
Broadly speaking, there are two approaches for decision making in uncertain environments: \emph{stochastic} optimization (SO) and \emph{robust} optimization (RO). SO relies on probabilistic models to explain the \se{uncertainties} in data and often results in solutions that are sensitive to these assumptions\footnote{We refer an interested reader to \cite{Birge2011} for a survey on stochastic modeling and techniques.}. 
On the other hand, RO incorporates a set-based deterministic model of the uncertainty such that the optimal solution protects against all realizations in the uncertainty set. 
Compared with SO, one significant advantage of RO is the \se{computational} tractability, which is important for the phase balancing problem to be applicable in real-world scenarios. Moreover, it has been observed that robust solutions are competitive with the deterministic solutions \se{in terms of cost, while being more robust to unplanned uncertainties in the data. Robust optimization also does not need to assume any probabilistic information about the uncertain quantities \cite{Bertsimas2004}.} 


We consider the following \se{row-uncertain} robust linear optimization problem, where the row vectors $\alpha_i$ are uncertain in each constraint:
\begin{subequations}
\label{opt:robust_linear_program}
\begin{align}
\minimize_{x} \quad & \gamma^\intercal x \\
\text{subject to} \quad & \alpha_i^\intercal x \le \beta_i, ~\forall \alpha_i \in \mathcal{U}_i, \\
& i = 1,2,\cdots, m. \nonumber
\end{align}
\end{subequations}
Formulation (\ref{opt:robust_linear_program}) seeks an optimal solution $x \in \mathbb{R}^n$ that is feasible to $m$ linear uncertain constraints $\alpha_i^\intercal x \le \beta_i$, in which the uncertain vector of parameters $\alpha_i$ can take any values from the uncertainty set $\mathcal{U}_i$. A common choice is the polyhedral uncertainty set defined as
\begin{equation}
\mathcal{U}_i := \{\alpha_i: H_i \alpha_i \le h_i \}, i = 1,2,\cdots,m,
\end{equation}
where $H_i \in \mathbb{R}^{k\times n}$ and $h_i \in \mathbb{R}^k$ depict $k$ inequalities that define a polyhedron. Such uncertainty sets have been successful in capturing insights from probability theory to obtain more realistic models. For instance, if the data is generated independently from a probability distribution then the well-known central limit theorem states that the appropriately normalized average of variables tends to a normal distribution. The central limit theorem can be written as a polyhedral uncertainty set that protects against all realizations of data that satisfy the central limit theorem \cite{Bandi2012}. Its parameters can be set such that if the data was generated via a given probability distribution, then the uncertainty set captures provably $95\%-99\%$ of possible scenarios. This provides a clean way to incorporate probabilistic information. We refer the reader to \cite{BenTal2009} for a more detailed survey on robust optimization techniques.


In the definition of a polyhedral uncertainty set $\mathcal{U} = \{\alpha \in \mathbb{R}^n | H \alpha \leq h\}$ where $H \in \mathbb{R}^{k \times n}$, the constraint $\alpha^\intercal x \le \beta ~\forall \alpha \in \mathcal{U}$ is equivalent to
\begin{subequations}
\label{opt:robust_linear_constraint}
\begin{align}
b \geq \maximize_{\alpha} \quad & x^\intercal \alpha & \\
\text{subject to}\quad & H \alpha \le h. \label{dualvar}
\end{align}  
\end{subequations}
Let $p \in \mathbb{R}_+^k$ be the dual variable for \eqref{dualvar}. Then the dual linear program of \eqref{opt:robust_linear_constraint} is:
\begin{subequations}
\label{opt:robust_linear_constraint_dual}
\begin{align}
\minimize_{p} \quad & h^\intercal p\\
\text{subject to} \quad& H^\intercal p = x,\\
&p \geq 0.
\end{align}
\end{subequations}
By weak duality, any feasible solution $p$ of (\ref{opt:robust_linear_constraint_dual}) for a given $x$ provides a lower bound to (\ref{opt:robust_linear_constraint}), i.e. $h^\intercal p \leq^{(*)} \max_{\alpha \in \mathcal{U}} \alpha^\intercal x \leq b$, and the inequality (*) is tight for the optimal solution of the dual formulation in (\ref{opt:robust_linear_constraint_dual}), by strong duality. 



Therefore the uncertain constraints $\alpha_i^\intercal x \le \beta_i~ \forall \alpha_i \in \mathcal{U}_i$ are equivalent to the following deterministic constraints:
\begin{eqnarray*}
  h_i^\intercal p_i \le \beta_i, \quad H_i^\intercal p_i = x,\quad p_i \ge 0, 
\end{eqnarray*}
where each $p_i \in \mathbb{R}_+^k$ is \se{a vector of auxiliary variables} corresponding to the $i$th constraint in (\ref{opt:robust_linear_program}). The robust formulation  (\ref{opt:robust_linear_program}) with polyhedral uncertainty sets $\mathcal{U}_i$ is then equivalent to the following linear program \cite{bertsimas_theory_2011}:
\begin{subequations}
\label{opt:robust_linear_program_eq_polyhedral}
\begin{align}
\minimize_x \quad & \gamma^\intercal x \\
\text{subject to}\quad & h_i^\intercal p_i \le \beta_i,  \\
& H_i^\intercal p_i = x, \\
& p_i \in \mathbb{R}_+^k, i = 1,2,\cdots m.
\end{align}
\end{subequations}

One major advantage of using a polyhedral uncertainty set is its computational tractability. The reformulation of the robust linear program (\ref{opt:robust_linear_program}) as the deterministic linear program (\ref{opt:robust_linear_program_eq_polyhedral}) involves a few more variables and this does not increase the \se{overall} computational complexity \cite{bertsimas_theory_2011}.

\section{Formulations of Phase Balancing Problems} 
\label{sec:phase_balancing}
\subsection{Nomenclature} 
Time dependent variables are represented with $\cdot[t]$, e.g. $d[t]$ is the demand at time $t$.
Matrices are represented using capital letters and uncertainty sets are in calligraphic font. $|\cdot|$ is the absolute value function and $^\intercal$ denotes transpose of matrices or vectors. By $\mathbf{1}$, we mean the vector of all ones in the appropriate dimension (typically $n$ in this paper, e.g. in \eqref{opt:det_phase_balancing_diff_a} $\mathbf{1} \in \mathbb{R}^n$).  

\subsection{Deterministic Phase Balancing} 
\label{sub:deterministic_phase_balancing}
We briefly review the conventional formulation of phase balancing in this subsection. Formulation (\ref{opt:det_phase_balancing}) presented below is a slight variation of the original one in \cite{zhu_phase_1998}.
\begin{subequations}
\label{opt:det_phase_balancing}
\begin{align}
\vspace{-0.2cm}
\minimize_{a,b, c, u_a, u_b, u_c} \quad & \max\{u_a, u_b, u_c\} \label{opt:det_phase_balancing_obj}\\
\text{subject to}\quad &  u_{a} =  | d^\intercal (a - \frac{\mathbf{1}}{3})|, \label{opt:det_phase_balancing_diff_a}\\
& u_{b} = | d^\intercal (b - \frac{\mathbf{1}}{3})|,  \label{opt:det_phase_balancing_diff_b}\\ 
& u_{c} = | d^\intercal (c - \frac{\mathbf{1}}{3})|, \label{opt:det_phase_balancing_diff_c}\\
& a + b + c = \mathbf{1},  \label{opt:det_phase_balancing_sum1}\\
& a, b, c \in \{0,1\}^n, \label{opt:det_phase_balancing_binary}\\
& u_{a} ,u_{b} , u_{c} \in \mathbb{R}_+.\label{opt:det_phase_balancing_pos_u}
\end{align}
\end{subequations}
Phase balancing aims at finding the most balanced assignment of $n$ loads $d \in \mathbb{R}^n$ to three phases (A,B,C).
Phase balancing commonly relies on phase swapping (or re-phasing) actions to reduce imbalances. Phase swapping typically happens at the feeder level, during maintenance or restoration periods \cite{zhu_phase_1998}. Phase swapping actions are depicted by decision variables $a$, $b$ and $c$, all of which are binary vectors \se{with dimension equal to the number of loads}, where $a_i=1$ (similarly, $b_i, c_i=1$) denotes load $d_i$ is assigned to phase A (similarly, to phase B, C), and $0$ indicates $d_i$ is not assigned to that phase.
Constraint (\ref{opt:det_phase_balancing_sum1}) ensures that each load must be assigned to exactly one phase\footnote{For simplicity, we only consider single-phase loads in this paper. Extensions to multi-phase loads are in Section \ref{sub:Multi_phase_loads}.}. Variables $u_a, u_b, u_c$ represent single-phase imbalances, namely the difference of load on phase A (B,C) from the uniformly balanced case $d^\intercal \mathbf{1}/3$. The objective (\ref{opt:det_phase_balancing_obj}) is to minimize the largest imbalance amongst the three phases. The original formulation in \cite{zhu_phase_1998} minimizes the largest differences between any two phases, i.e. 
\begin{align}
\label{opt:original_phase_balancing}
\minimize_{a, b, c \in \{0,1\}^n} \quad \max\{|d^\intercal (b-a)|, |d^\intercal (b-c)|, |d^\intercal (a-c)|\}.
\end{align}
These two formulations are closely related in the following sense. Let the total loads assigned to phases A, B and C be $x, y$ and $z$ respectively, and let the total overall load be $x+y+z=d^\intercal \mathbf{1}$. Without loss of generality, let $x\leq y\leq z$ which implies $x \leq d^\intercal \mathbf{1}/3 \leq z$. So, the objective value of \eqref{opt:original_phase_balancing} will be $z-x$ whereas the objective value of (\ref{opt:det_phase_balancing}) for such an assignment will be $\max\{d^\intercal \mathbf{1}/3 - x, z - d^\intercal \mathbf{1}/3\}$, but note that $z - x \leq 2 \max \{\mathbf{1}/3 - x, z - d^\intercal \mathbf{1}/3\}$. Therefore, the optimal solution of \eqref{opt:original_phase_balancing} will be at most twice the optimal solution of \eqref{opt:det_phase_balancing} (and similarly, optimal solution of \eqref{opt:original_phase_balancing} is at least the optimal solution of \eqref{opt:det_phase_balancing}). Further, we believe that our formulation in \eqref{opt:det_phase_balancing} meets the intuitive notion of phase balancing better than \eqref{opt:original_phase_balancing}. To see this, consider a total given demand of 21 kW. Formulation \eqref{opt:original_phase_balancing} does not differentiate between the assignments 2,9,10 kW and 3,7,11 kW on each phase. For either assignment, the maximum difference between the assigned loads is 8 kW. However, our formulation \eqref{opt:det_phase_balancing} would prefer 3,7,11 as a solution since it minimizes the maximum deviation from the average.
The absolute value constraints \eqref{opt:det_phase_balancing_diff_a}, \eqref{opt:det_phase_balancing_diff_b} and \eqref{opt:det_phase_balancing_diff_c} can be reformulated to obtain an equivalent mixed integer linear program \cite{ben-tal_robust_2009}:
\vspace{-0.2cm}
\begin{subequations}
\label{opt:det_eq_phase_balancing}
\begin{align}
\minimize_{a,b,c,u} \quad & u \label{opt:det_eq_phase_balancing_obj}\\
\text{subject to}\quad &  -u \le d^\intercal (a - \frac{\mathbf{1}}{3}) \le u, \label{opt:det_eq_phase_balancing_diff_a}\\
& -u \le d^\intercal (b - \frac{\mathbf{1}}{3}) \le u, \label{opt:det_eq_phase_balancing_diff_b}\\ 
& -u \le d^\intercal (c - \frac{\mathbf{1}}{3}) \le u, \label{opt:det_eq_phase_balancing_diff_c}\\
& a + b + c = \mathbf{1}, \\
& a, b, c \in \{0,1\}^n, u \in \mathbb{R}_+. 
\end{align}
\end{subequations}

\subsection{Robust Phase Balancing} 
\label{sub:robust_phase_balancing}
In deterministic phase balancing problem (\ref{opt:det_phase_balancing}), load vector $d$ represents the average load level during a long period, without any uncertainties. 
Motivated by the rapid growth of highly variable resources in distribution systems, we connect conventional phase balancing with robust optimization and formulate the following robust phase balancing problem: 
\begin{subequations}
\label{opt:robust_phase_balancing}
\begin{align}
\minimize_{u, a, b, c} \quad & u \label{opt:robust_eq_phase_balancing_obj}\\
\text{subject to}\quad &  -u \le d^\intercal (a - \frac{\mathbf{1}}{3}) \le u, \forall d \in \mathcal{D}, \label{opt:robust_phase_balancing_diff_a}\\
& -u \le d^\intercal (b - \frac{\mathbf{1}}{3}) \le u, \forall d \in \mathcal{D},  \label{opt:robust_phase_balancing_diff_b}\\ 
& -u \le d^\intercal (c - \frac{\mathbf{1}}{3}) \le u, \forall d \in \mathcal{D},  \label{opt:robust_phase_balancing_diff_c}\\
& a + b + c = \mathbf{1}, \\
& a, b, c \in \{0,1\}^n, u \in \mathbb{R}_+.
\end{align}
\end{subequations}
The major difference \se{between} robust phase balancing (\ref{opt:robust_phase_balancing}) \se{and} the deterministic version (\ref{opt:det_eq_phase_balancing}) is that:
instead of seeking solutions $(a,b,c)$ that are feasible for the average or expected load vector $d$, (\ref{opt:robust_phase_balancing}) seeks solutions robust to all realizations of $d$ in an uncertainty set $\mathcal{D}$. The uncertainty set $\mathcal{D}$ can be constructed using historical data or approximated with prior knowledge. 
 
Similar to Section \ref{sec:robust_optimization},
formulation (\ref{opt:robust_phase_balancing}) with polyhedral uncertainty set $\mathcal{D} = \{d: H d \le h \}$ can be rewritten as an MILP (\ref{opt:robust_phase_balancing_eq_polyhedral}).
\begin{subequations}
\label{opt:robust_phase_balancing_eq_polyhedral}
\begin{align}
\minimize_{p, q, a, b, c, u} \quad & u \label{opt:robust_phase_balancing_obj}\\
\text{subject to}\quad &  h^\intercal p_a \le u, \quad H^\intercal p_a = a - \frac{\mathbf{1}}{3}, \\
&  h^\intercal q_a \le u, \quad H^\intercal q_a = \frac{\mathbf{1}}{3} -a, \\
& h^\intercal p_b \le u, \quad H^\intercal p_b = b - \frac{\mathbf{1}}{3}, \\
& h^\intercal q_b \le u, \quad H^\intercal q_b = \frac{\mathbf{1}}{3} -b, \\
& h^\intercal p_c \le u, \quad H^\intercal p_c = c - \frac{\mathbf{1}}{3}, \\
& h^\intercal q_c \le u, \quad H^\intercal q_c = \frac{\mathbf{1}}{3} -c, \\
& a + b + c = \mathbf{1}, \\
& a, b, c \in \{0,1\}^n, u \in \mathbb{R}_+, \\
& p_a, p_b, p_c, q_a, q_b, q_c \in \mathbb{R}_+^k. 
\end{align}
\end{subequations}
where $p_a,p_b,p_c$ and $q_a,q_b,q_c$ are auxiliary variables.



\subsection{Robust Look-ahead Phase Balancing} 
\label{sub:robust_look_ahead_phase_balancing} 
The problem formulated in Section \ref{sub:robust_phase_balancing} considers only a single snapshot (e.g. one hour) decision making for robust phase balancing. However, one key component of costs comes from frequent phase swapping actions of loads. Therefore, it \se{is} important to consider the phase balancing problem in a multi-time-horizon setting. We formulate it as a robust look-ahead phase balancing problem, much like the usual practices in \cite{thatte_robust_2014,xie_short-term_2014}.

In the following formulation \se{\eqref{opt:robust_dynamic_phase_balancing}}, we consider a two-period moving horizon phase balancing decision making. For example, each snapshot (i.e., $t$) could signify two hours in the day, a period $1, \hdots, T_1$ (i.e. a day) consists of 12 snapshots, and the moving horizon might consist of two days. The objective function (i.e., phase imbalances) is defined over the two periods combined. However, \emph{the decisions are only implemented for the first period.} The reason for doing so is justified by the engineering insight that information gets more accurate as we get closer to real-time operations. Therefore, the decision made for period two is only advisory but not implemented. 

\begin{figure}[htbp]
  \centering
  \includegraphics[width=\linewidth]{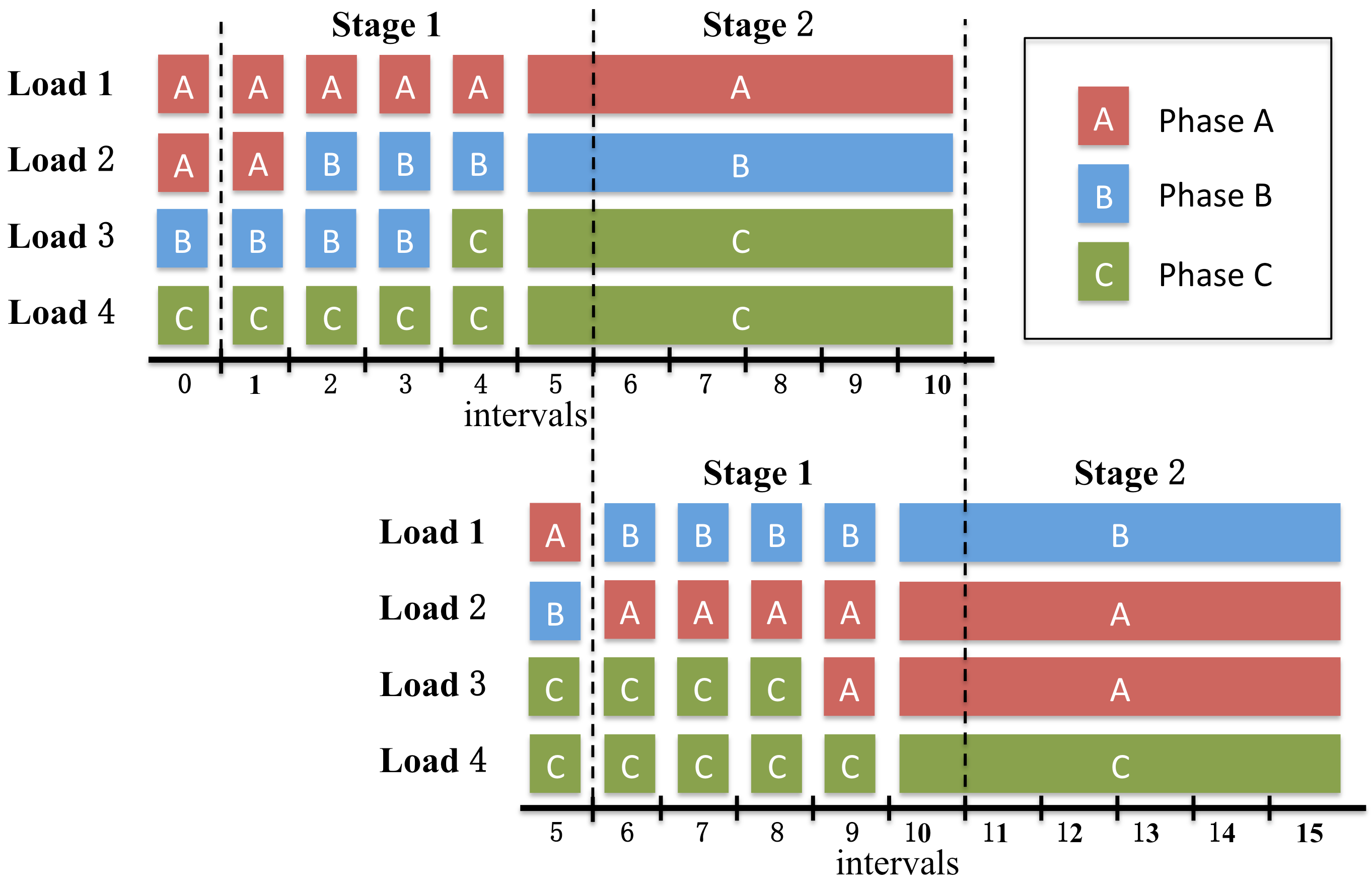}
  \caption{Illustration of the Look-ahead Operation Framework (every block represents the phase assignment of a load at each snapshot)}
  \label{fig:fig_rDPB}
\end{figure}

An illustrative example with 4 loads over 10 intervals is provided in Figure \ref{fig:fig_rDPB}. 
A robust look-ahead phase balancing problem is solved for \se{intervals} 1 to 10 given initial load assignments at interval 0. \se{Two phase swaps are implemented in the illustrated solution (top of Figure \ref{fig:fig_rDPB}:} load 2 is swapped to phase B at the beginning of interval 2, load 3 is swapped to phase C at the beginning of interval 4. \se{At the end of interval 5,} uncertainty sets for interval 6 to 15 \se{are} constructed using updated load forecast and another robust look-ahead phase balancing problem is solved \se{for intervals 6 to 10}. \se{The solution implemented performs three phase swaps in intervals 6 to 10, as shown at the bottom of Figure \ref{fig:fig_rDPB}.}

\se{We next provide a mixed integer formulation to solve the two-period moving horizon phase balancing problem.}

\begin{small}
\begin{subequations}
\allowdisplaybreaks
\label{opt:robust_dynamic_phase_balancing}
\begin{align}
\minimize \quad & u + \lambda v \\
\text{subject to}\quad &  -u \le (d[t])^\intercal (a[t] - \frac{\mathbf{1}}{3}) \le u, \forall d[t] \in \mathcal{D}_t, \label{opt:robust_dynamic_phase_balancing_diff_a1}\\
&  -u \le (d[t])^\intercal (b[t] - \frac{\mathbf{1}}{3}) \le u, \forall d[t] \in \mathcal{D}_t, \label{opt:robust_dynamic_phase_balancing_diff_b1}\\
&  -u \le (d[t])^\intercal (c[t] - \frac{\mathbf{1}}{3}) \le u, \forall d[t] \in \mathcal{D}_t, \label{opt:robust_dynamic_phase_balancing_diff_c1}\\
& \qquad \qquad \qquad \qquad \qquad \quad t = 1,2,\cdots, T_1 \nonumber \\
&  -v \le (d[t])^\intercal (a[T_1] - \frac{\mathbf{1}}{3}) \le v, \forall d[t] \in \mathcal{D}_t, \label{opt:robust_dynamic_phase_balancing_diff_a2}\\
&  -v \le (d[t])^\intercal (b[T_1] - \frac{\mathbf{1}}{3}) \le v, \forall d[t] \in \mathcal{D}_t, \label{opt:robust_dynamic_phase_balancing_diff_b2}\\
&  -v \le (d[t])^\intercal (c[T_1] - \frac{\mathbf{1}}{3}) \le v, \forall d[t] \in \mathcal{D}_t, \label{opt:robust_dynamic_phase_balancing_diff_c2}\\
& \qquad \qquad \qquad \quad t = T_1+1,T_1+2,\cdots, T_2 \nonumber \\
& \sum_{t=1}^{T_1} \Big( \mathbf{1}^\intercal |a[t] - a[t-1]| + \mathbf{1}^\intercal |b[t] - b[t-1]| \nonumber\\
& \qquad \qquad \qquad \quad + \mathbf{1}^\intercal |c[t] - c[t-1]| \Big) \le 2s, \label{opt:robust_dynamic_phase_balancing_switchings} \\
& a[t] + b[t] + c[t] = \mathbf{1},  \\
& a[t], b[t], c[t] \in \{0,1\}^n, u,v \in \mathbb{R}_+, \\
& \qquad \qquad \qquad \qquad \qquad \quad t = 1,2,\cdots,T_1. \nonumber
\end{align}
\end{subequations}
\end{small}

\se{In the above formulation \eqref{opt:robust_dynamic_phase_balancing}}, the first period consists of $T_1$ snapshots ($t=1,2,\cdots,T_1$). It determines the phase swapping actions to be implemented. Similar to previous formulations, $a_i[t]=1$ indicates load $d_i[t]$ is assigned to phase A at time $t$ ($t=1,2,\cdots,T_1$). (\ref{opt:robust_dynamic_phase_balancing_diff_a1})-(\ref{opt:robust_dynamic_phase_balancing_diff_c1}) are robust constraints for period 1. It is worth noting that each snapshot has its own uncertainty set $d[t] \in \mathcal{D}_t$.
This allows (\ref{opt:robust_dynamic_phase_balancing}) to take advantage of the temporal patterns of uncertain loads.
As illustrated in Figure \ref{fig:fig_rDPB}, no phase swapping actions are considered for period 2.
Formulation (\ref{opt:robust_dynamic_phase_balancing}) seeks \emph{fixed} load assignments with small phase imbalances for period 2. The decision variables of the second period are $a[T_1]$, $b[T_1]$ and $c[T_1]$.
Constraints (\ref{opt:robust_dynamic_phase_balancing_diff_a2})-(\ref{opt:robust_dynamic_phase_balancing_diff_c2}) \se{relate to decisions in} period 2.

\se{We do not allow} phase swapping actions \se{in} the second period of (\ref{opt:robust_dynamic_phase_balancing}) for two \se{important} reasons: (a) uncertainties for the second period could be significantly larger than in the first one, over-optimization with large uncertainties might lead to conservative solutions; (b) the problem size will be twice larger if \se{we consider} phase swapping in both periods \se{thus hurting performance}. Recall that phase balancing is an MILP, the computational burden could be prohibitive\footnote{We actually tested the case in which phase swapping is considered in both periods. Gurobi \cite{gurobi_optimization_gurobi_2016} took 12 hours to converge and the solution was \se{comparable to the current formulation in} (\ref{opt:robust_dynamic_phase_balancing}).}

Variables $u$ and $v$ denote the largest single-phase imbalance that occurs in the two periods, respectively. 
Choosing a proper value of parameter $\lambda \in \mathbb{R}_+$ could achieve a balance between the optimality \se{in short term and long term}.

Given current industrial practice, swapping loads from one phase to another typically requires manual operations, which incurs extra costs on human resources, maintenance expenses and planned outage duration \cite{zhu_phase_1998}. Constraint (\ref{opt:robust_dynamic_phase_balancing_switchings}) limits the maximum number of phase swapping actions in the first period. Parameter $s$ denotes the budget of swapping actions. Without constraint (\ref{opt:robust_dynamic_phase_balancing_switchings}), a large amount of phase swapping actions could be recommended, which is not affordable for utility companies \cite{zhu_phase_1998}.

For each snapshot $t=1,2,\cdots, T_2$, the polyhedral uncertainty set is defined as
\begin{equation}
  \mathcal{D}_t = \{d[t]: H_t d[t] \le h_t \}
\end{equation}

By introducing auxiliary variables, (\ref{opt:robust_dynamic_phase_balancing}) is equivalent to an MILP (\ref{opt:robust_dynamic_phase_balancing_eq}).

It is worth mentioning that a recent paper \cite{sun_phase_2016} proposes a related but different approach with stochastic optimization. It minimizes the expected loss function over a time horizon with respect to uncertainties from loads and electricity prices. While its decision variables \se{denote the} charging and discharging rates of energy storages, load assignments remain unchanged and no phase swapping actions are considered.

\section{Case Study} 
\label{sec:case_study}
\subsection{Load Data} 
\label{sub:load_data}
The load profiles are from dataset ``R1-12.47-4'' of \cite{hoke_steady-state_2013}. It models a heavily populated suburban area composed mainly of single family homes and heavy commercial loads \cite{schneider_modern_2008}. The dataset ``R1-12.47-4'' is populated with hourly averaged load data from a utility company in the West Coast of the United States \cite{hoke_steady-state_2013}. The original dataset is publicly available on \url{catalog.data.gov}.
The dataset contains $74$ hourly load profiles of 365 days. We use the first 30 days and scale them randomly to avoid identical load profiles.

\begin{figure}[h]
	\centering
  \subfloat[average daily load profiles with standard deviations (each color represents one load)]{\label{fig:R1-1247-4-load-profile-all} \includegraphics[width= \linewidth]{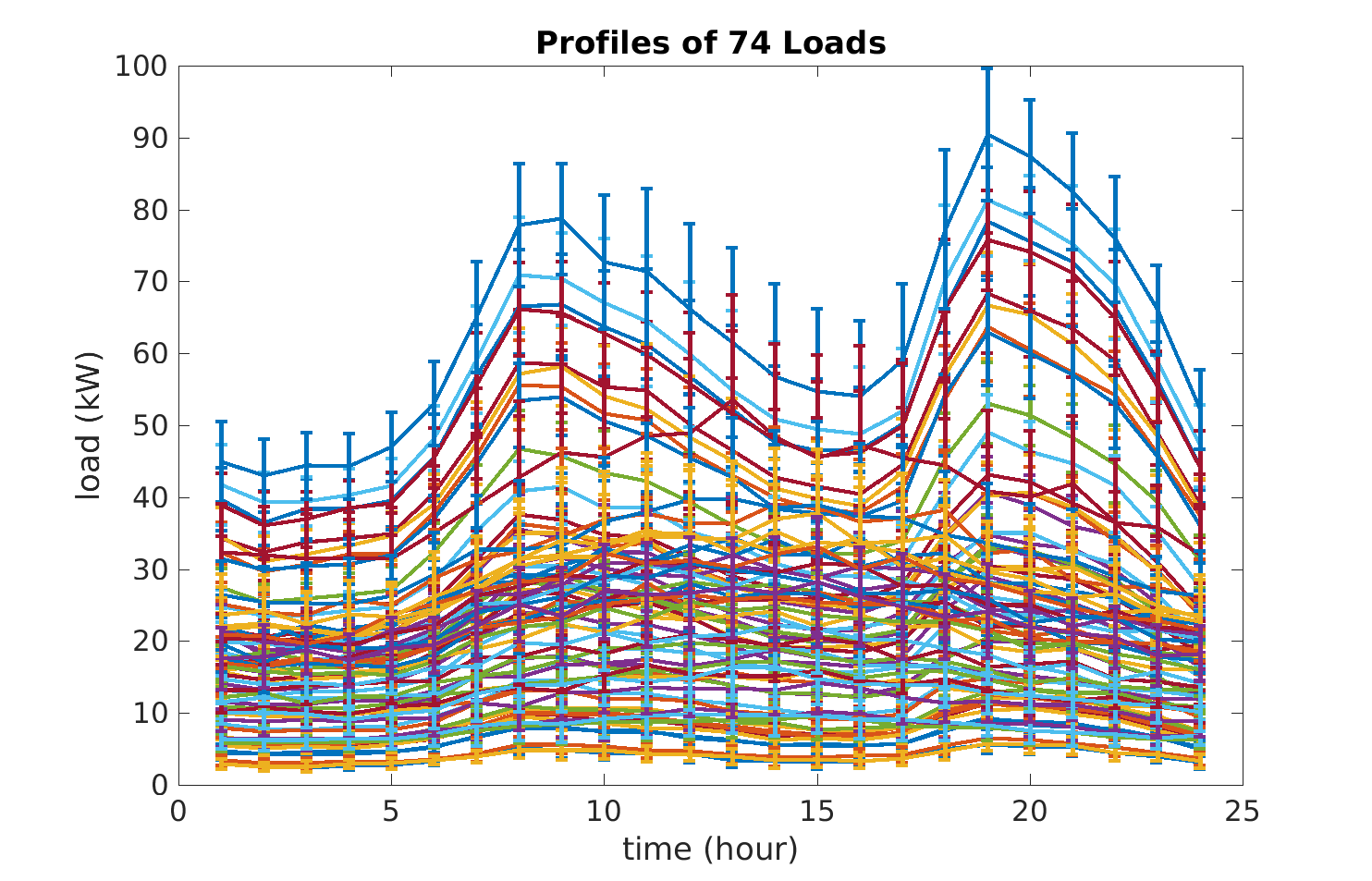}} \\
  \subfloat[profiles of load 16 (different colors represent different days)]{\label{fig:R1-1247-4-load-profile-one} \includegraphics[width= \linewidth]{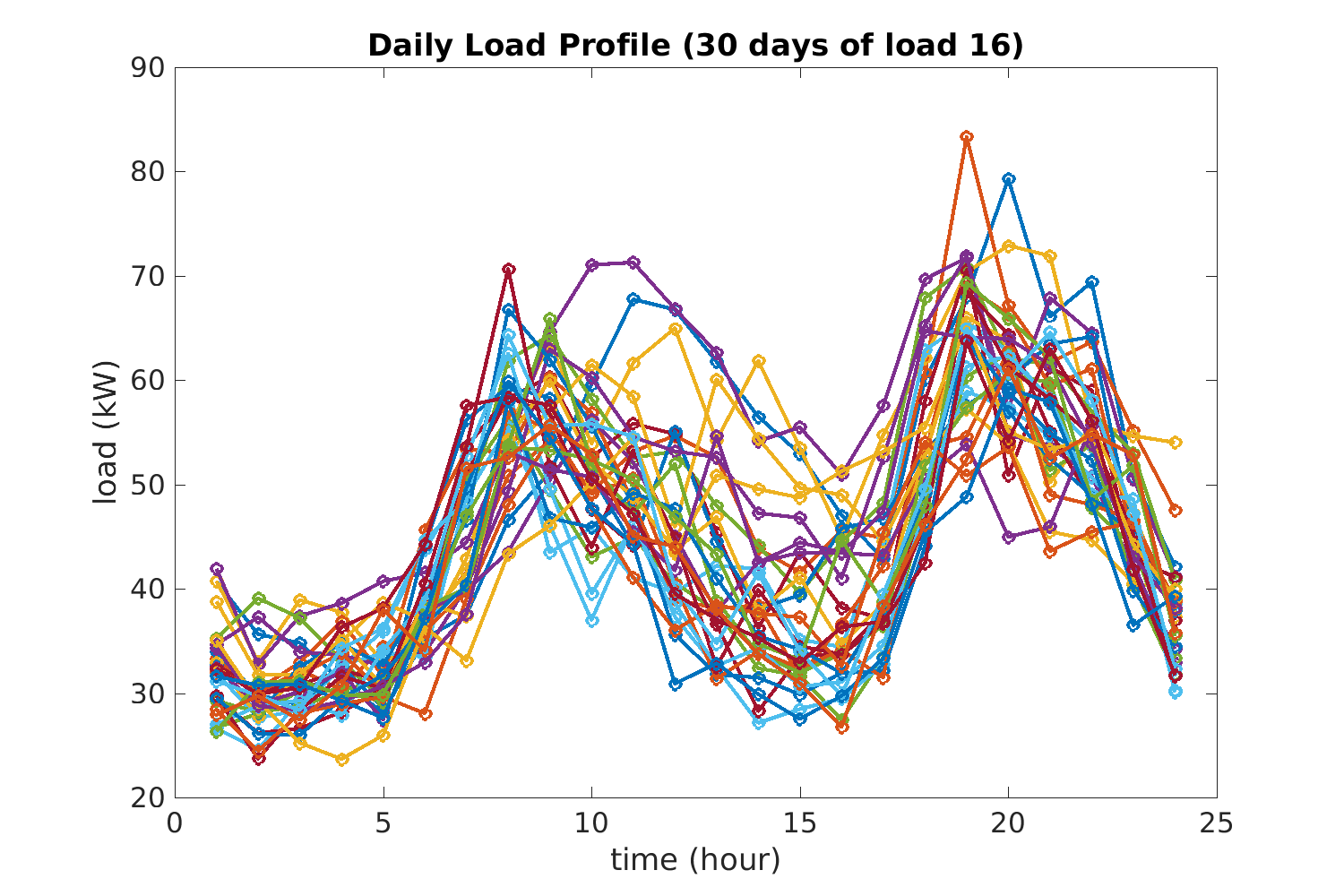}}    
	\caption{Modified Load Dataset ``R1-1247-4''}
	\label{fig:R1-1247-4-load-profile}
\end{figure}

Figure \ref{fig:R1-1247-4-load-profile} visualizes the modified dataset.

\subsection{Construct Uncertainty Set} 
\label{sub:construct_uncertainty_set}
In order to demonstrate the benefits of robustification, we use the following polyhedral uncertainty sets for the \se{robust Phase Balancing (r-PB) and robust Look-ahead Phase Balancing (r-LAPB) problems:}
\begin{eqnarray}
\label{eqn:def_uncertainty_set_rSPB}
  \mathcal{D} = \{d \in \mathbb{R}^n: \hat{d} \le d - \overline{d} \le \hat{d} \}
\end{eqnarray}	
\begin{equation}
\label{eqn:def_uncertainty_set_rDPB}
	\mathcal{D}_t = \{d[t] \in \mathbb{R}^n:(1-\rho_t) \overline{d[t]} \le d[t] \le (1+\rho_t) \overline{d[t]}   \}	
\end{equation}
where $\overline{d} \in \mathbb{R}^n$ or $\overline{d[t]}\in \mathbb{R}^n$ represent the average load or forecast value, and $\hat{d}\in \mathbb{R}^n$ denotes the largest deviation of load $d$. Problem r-PB  (\ref{opt:robust_phase_balancing}) with $\hat{d}=0$ is equivalent with \se{deterministic Phase Balancing (d-PB)} (\ref{opt:det_phase_balancing}). Values of $\overline{d}$, $\overline{d}[t]$ and $\hat{d}$ are estimated from the modified ``R1-1247-4'' dataset. These uncertainty sets can be viewed as simple relaxations of the central limit theorem based sets (which can risk the solution being too conservative), but they already show a significant improvement in our experiments compared to deterministic solutions.

For r-LAPB, the level of robustness $\rho_t$ depends on the forecast accuracy or confidence. Larger $\rho_t$ indicates lower forecast accuracy. Definition of $\mathcal{D}_t$ in (\ref{eqn:def_uncertainty_set_rSPB})-(\ref{eqn:def_uncertainty_set_rDPB}) assumes that the load forecast is unbiased and bounded by $\rho_t$. For r-LAPB, $\rho_t$ in the first period (i.e. 24 hours) is set to be $10\%$ ($t=1,2,\cdots,24$) and $\rho_t = 30\%$ for the second period ($t=25,26,\cdots,48$).

For d-PB (\ref{opt:det_phase_balancing}), load vector $d$ is the average hourly load of 30 days. There is no uncertainty associated.

\begin{table}[tb]
	\caption{Parameters}
	\label{tab:parameters}
	\centering

	\begin{tabular}{cccccc}
	\hline

	\hline
	$n$ & $T_1$ & $T_2$&  $\lambda$ & $\rho_t$ (period1) & $\rho_t$ (period2)\\
	\hline
	74 & 24   & 48    & 1/3 & 10\% & 30\% \\
	\hline

	\hline
	\end{tabular}
\end{table}

\subsection{Simulation Results} 
\label{sub:simulation_results}
\begin{figure}[htbp]
  \centering
  \includegraphics[width=\linewidth]{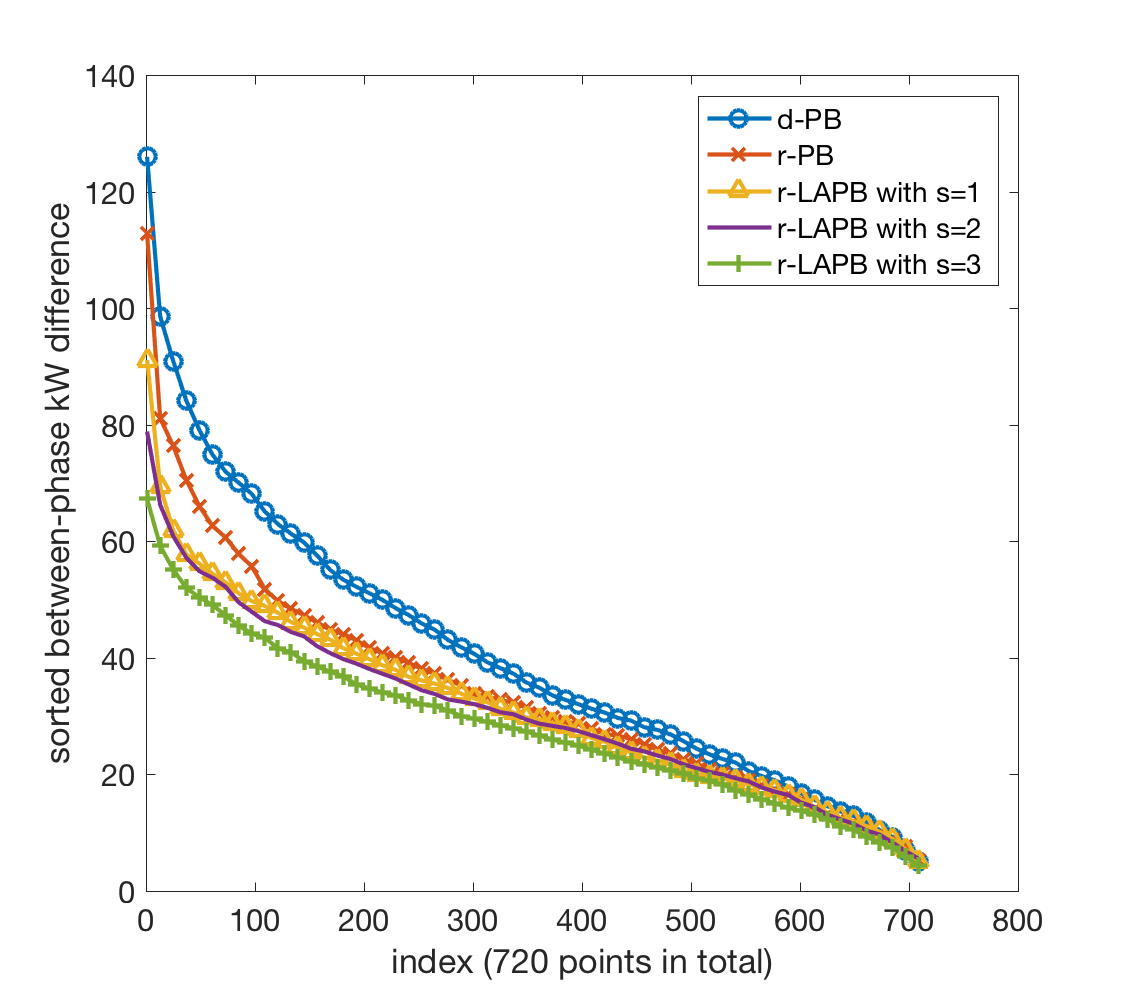}
  \caption{Sorted between-phase kW differences}
  \label{fig:R1-1247-4-sorted-max-kW-diff}
\end{figure}

Simulations are performed on a desktop with Intel i7-2600 8-core CPU@3.40GHz and 16GB memory. The phase balancing problems are solved using YALMIP \cite{lofberg_automatic_2012,lofberg_yalmip_2004} and Gurobi \cite{gurobi_optimization_gurobi_2016}. The optimality gap of every solution is smaller than $0.1\%$.
Key results are presented in Figure \ref{fig:R1-1247-4-sorted-max-kW-diff} and Table \ref{tab:key_results}.

The performance of three formulations are evaluated using three metrics: between-phase kW difference $\omega$, single-phase kW difference $\nu$ and single-phase percentage difference $\upsilon$, which are defined below:

\begin{small}
\begin{flalign*}
  \omega &:= \max\{|d^\intercal(a-b)|,|d^\intercal(a-c)|,|d^\intercal(b-c)|\}, \\
  \nu &:= \max\{|d^\intercal(a- \mathbf{1}/3)|,|d^\intercal(b-\mathbf{1}/3)|,|d^\intercal(c-\mathbf{1}/3)|\}, \\
  \upsilon &:= \max\{|1- \frac{3d^\intercal a}{d^\intercal \mathbf{1}}|,|1- \frac{3d^\intercal b}{d^\intercal \mathbf{1}}|,|1- \frac{3d^\intercal c}{d^\intercal \mathbf{1}}|\}. 
  \end{flalign*}
\end{small}

\begin{table*}[t]
  \caption{Comparison of Solutions}
  \label{tab:key_results}
  \centering
  \begin{tabular}{l|cccccccccccc}
  \hline

  \hline
                  & \multicolumn{3}{c}{between-phase (kW)}  & \multicolumn{3}{c}{single-phase (kW)} & \multicolumn{3}{c}{single-phase (\%)} & \multicolumn{3}{c}{runtime (s)}\\
  \textbf{Method} & \textbf{max} & \textbf{avg} & \textbf{std} & \textbf{max} & \textbf{avg} & \textbf{std}  & \textbf{max} & \textbf{avg} & \textbf{std} &  \textbf{max} & \textbf{avg} & \textbf{std} \\
  \hline
  d-PB           & 125.97 & 39.71 & 23.45 & 71.63 & 22.91 & 13.36 & 10.93 & 4.03 & 2.13 & - & 15.8 & - \\
  r-PB           & 112.79 & 33.55 & 18.78 & 63.66 & 19.27 & 10.87 & 10.56 & 3.44 & 1.87 & - & 1.3 & - \\
  r-LAPB {\tiny($s=1$)} & 91.09  & 31.12 & 15.81 & 49.19 & 17.91 & 9.11 & 10.62 & 3.25 & 1.69 & 109.0 & 56.2 & 18.3 \\
  r-LAPB {\tiny($s=2$)} & 78.83  & 30.39 & 15.02 & 45.84 & 17.45 & 8.60 & 9.71  & 3.20 & 1.69 & 1317.0 & 285.2 & 246.7 \\
  r-LAPB {\tiny($s=3$)} & 67.28  & 27.77 & 13.74 & 40.69 & 16.00 & 7.93 & 8.59  & 2.97 & 1.65 & 7756.6 & 2380.5 & 1321.2 \\
  \hline

  \hline
  \end{tabular}
\end{table*}

Compared with d-PB, robust phase balancing (r-PB) reduces both between-phase and single-phase imbalances by around 11\%, the standard deviations of imbalances are reduced by more than 20\%. It is also worth mentioning that the \se{time} to solve r-PB is significantly reduced due to more restricted search space since the robust solutions must be feasible for all demand realizations. 

Table \ref{tab:key_results} shows that the imbalances could be significantly reduced by incorporating look-ahead operations. For example, r-LAPB with 3 swapping actions per day reduces both between-phase and single-phase kW differences by 30\% on average. 

Figure \ref{fig:R1-1247-4-sorted-max-kW-diff} and Table \ref{tab:key_results} also demonstrate the trade-off between performance and computation complexity. In general, more frequent phase swapping operations lead to less imbalances among phases, while the time of solving r-LAPB grows exponentially\footnote{The r-LAPB with $s=4$ typically requires around $10\sim 12$ hours to solve one instance.}. 
Figure \ref{fig:R1-1247-4-sorted-max-kW-diff} clearly shows the major improvement of performance happens at the stage of applying r-PB and r-LAPB with one swapping per day. Improvements of allowing more swapping actions are marginal at the cost of higher \se{computational} burden and possible extra \se{cost} on human resources and maintenance.

\begin{figure}[htbp]
  \centering
  \subfloat[r-PB]{\includegraphics[width= \linewidth]{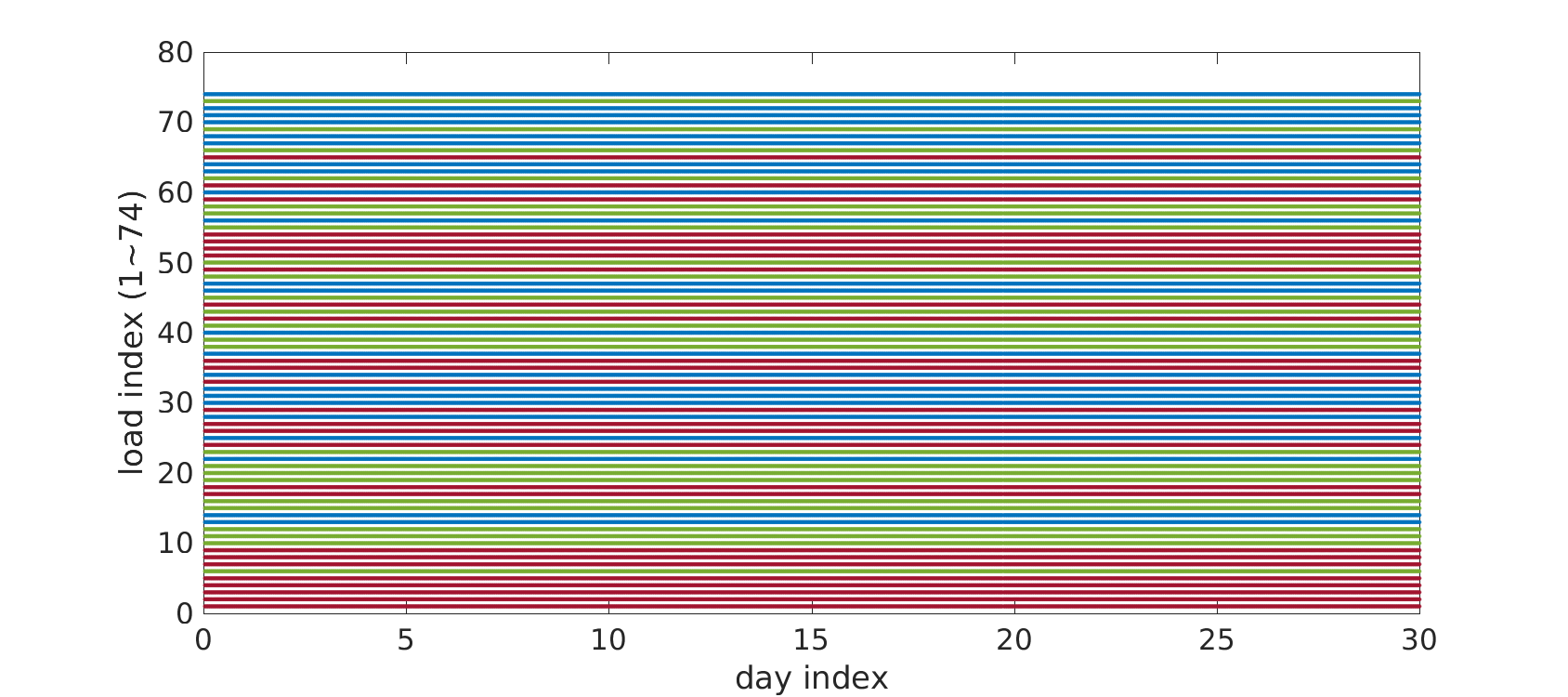}} \\
  \subfloat[r-LAPB ($s=1$)]{\label{R1-1247-4-rLAPB-box-rDPB-S=1-solution-phaseABC}\includegraphics[width= \linewidth]{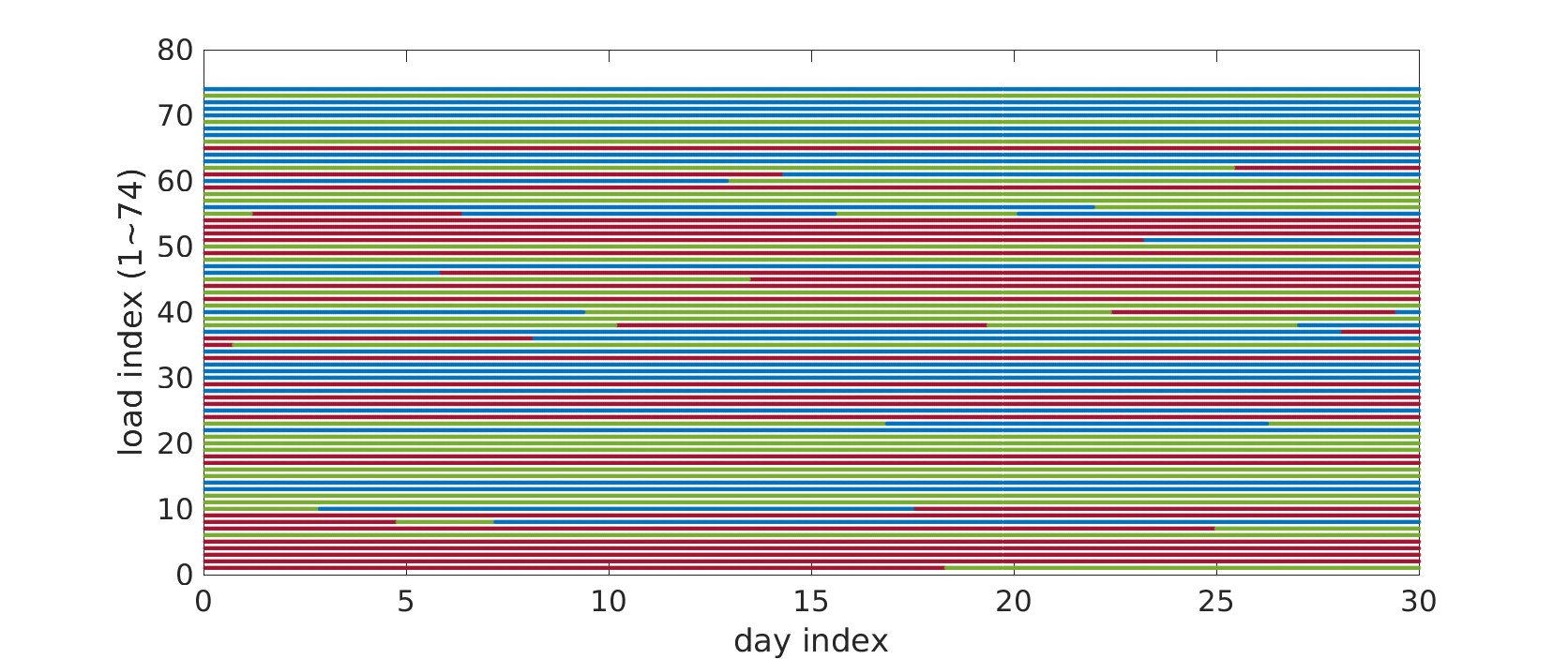}} \\
  \subfloat[r-LAPB ($s=2$)]{\label{R1-1247-4-rLAPB-box-rDPB-S=2-solution-phaseABC}\includegraphics[width= \linewidth]{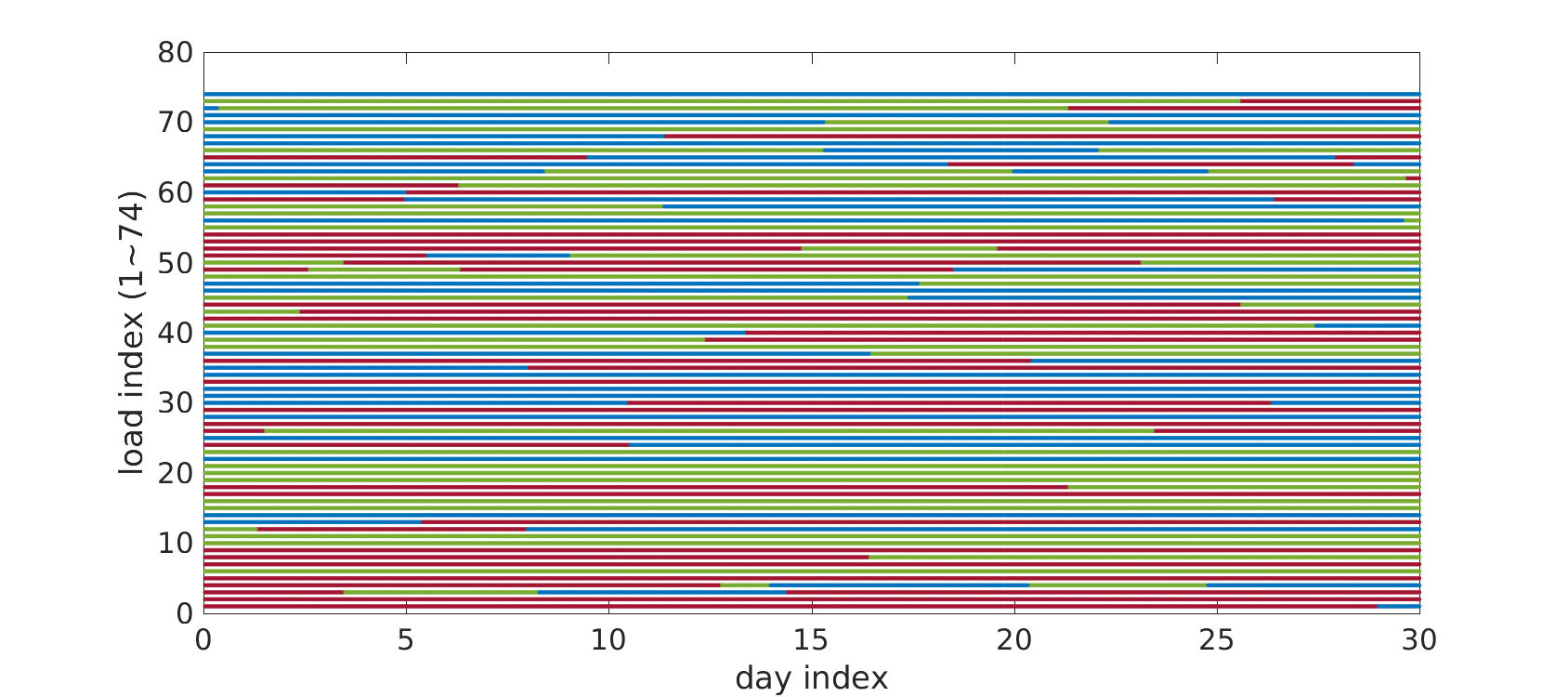}
  } \\
  \subfloat[r-LAPB ($s=3$)]{\label{R1-1247-4-rLAPB-box-rDPB-S=3-solution-phaseABC}\includegraphics[width= \linewidth]{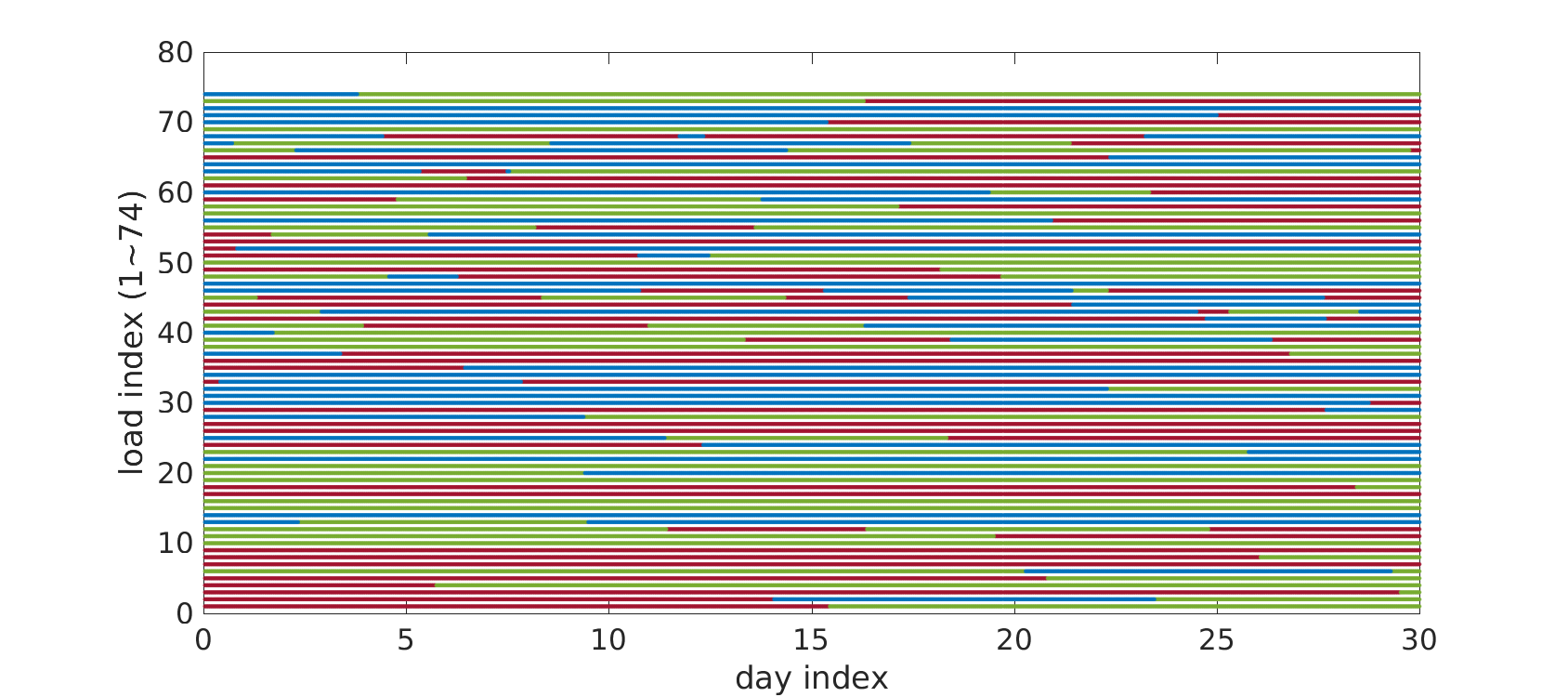}} 
  \caption{Display of Optimal Solutions (ABC phases are color-coded, red:phase A, blue:phase B, green:phase C)}
  \label{fig:R1-1247-4-rLAPB-box-rDPB-solution-phaseABC}
\end{figure}

\begin{figure}[htbp]
  \centering
  \includegraphics[width=\linewidth]{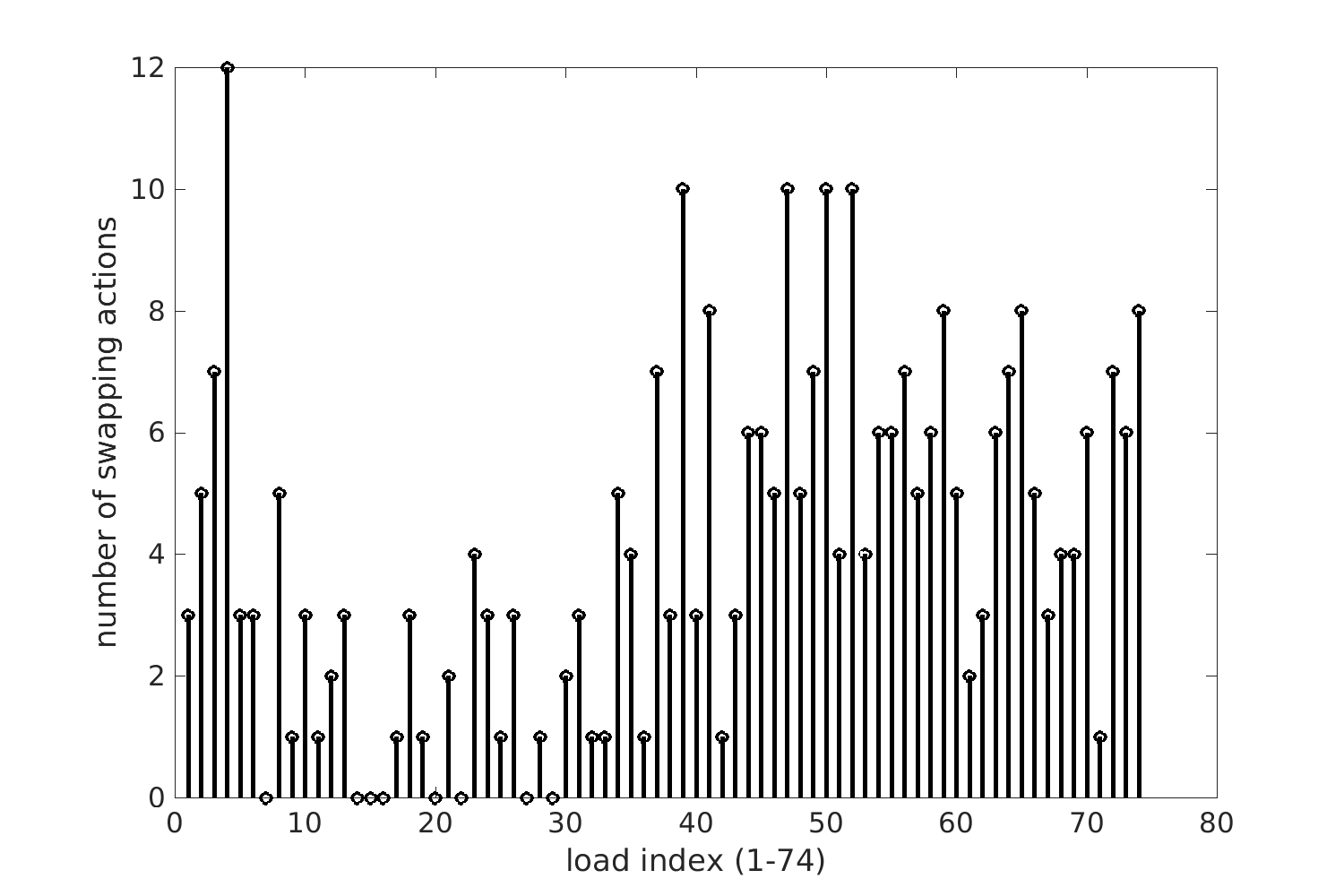}
  \caption{Phase Swapping Actions of Each Load (results of r-LAPB ($s=2$) running for 150 days)}
  \label{fig:R1-1247-4-rLAPB-num_switch-S=2}
\end{figure}

We also examine the optimal solution of r-PB and r-LAPB (Figure \ref{fig:R1-1247-4-rLAPB-box-rDPB-solution-phaseABC}). When allowing one swapping per day, all 28 switchings in 30 days happen on 18 out of 74 loads (Figure \ref{R1-1247-4-rLAPB-box-rDPB-S=1-solution-phaseABC}). Many loads remain unchanged and some loads have more frequent phase swapping operations than others. Figure \ref{fig:R1-1247-4-rLAPB-num_switch-S=2} demonstrates the case where two swapping actions are allowed per day, fives loads are swapped much more frequently than the others (17.3\% of 299 actions in 150 days, whereas the remaining 69 loads are switched only 3-4 times in 150 days on average.). Automatic phase swapping devices could be installed at these locations for more efficient and frequent responses.



\section{Discussions} 
\label{sec:discussions}
\subsection{Uncertainty Sets} 
\label{sub:uncertainty_sets}

In this paper, the uncertainty sets (\ref{eqn:def_uncertainty_set_rSPB})-(\ref{eqn:def_uncertainty_set_rDPB}) we use are a special case of polyhedral uncertainty sets. We do not capture yet potential correlations among different loads, as shown in Figure \ref{fig:R1-1247-4-load-profile-all}. 
Other choices of uncertainty sets might outperform current \se{ones} and reduce conservativeness, e.g. central limit theorem based polyhedral sets \cite{Bandi2012}, ellipsoidal uncertainty sets \cite{ben-tal_robust_1999}, cardinality constrained uncertainty sets \cite{bertsimas_robust_2003}, and constructing polyhedral uncertainty sets from data \cite{bertsimas_data-driven_2018}.

\subsection{Approximation Algorithms} 
\label{sub:approximation_algorithms}
All our formulations of phase balancing problems are mixed integer programs, which \se{are} in general computationally intractable. One of the classical problems in combinatorial optimization is {\it minimum makespan scheduling} that attempts to run a given set of jobs on a fixed number of parallel machines such that total time, i.e. the makespan, to complete jobs on any machine is minimized \cite{Lenstra1990}. Minimizing the maximum total load on any phase can then be viewed as makespan scheduling where the given set of jobs is simply the various loads, and the three parallel identical machines are the three phase lines. 
It is an open question to adapt known approximation algorithms for the minimum makespan scheduling problem (or to develop new methods) to the robust framework while incorporating switching costs. 
Deterministic phase balancing (d-PB) can also be seen as the optimization version of the $k$-partition problem \cite{sahni_p-complete_1976}, that attempts to divide $n$ integers into $k$ subsets such that the total sum of each subset is close to each other. This problem is a generalization of the three phase balancing problem, and might provide useful insights \se{as well}. 

\subsection{Multi-phase Loads} 
\label{sub:Multi_phase_loads}
It is easy to extend current phase balancing \se{problems} for the consideration of multi-phase loads. For deterministic phase balancing (\ref{opt:det_phase_balancing}), we could define variable $a^{(1)}$, $b^{(1)}$ and $c^{(1)}$ for single phase loads, $a^{(2)}$, $b^{(2)}$ and $c^{(2)}$ for two-phase loads, $a^{(3)}$, $b^{(3)}$ and $c^{(3)}$ for loads connecting to all three phases. Instead of constraint (\ref{opt:det_phase_balancing_sum1}), we have the following constraints:
\begin{eqnarray*}
	a^{(1)} + b^{(1)} + c^{(1)} = \mathbf{1},& \\
	a^{(2)} + b^{(2)} + c^{(2)} = 2\cdot\mathbf{1},& \text{ and }\\
	a^{(3)} + b^{(3)} + c^{(3)} = 3\cdot\mathbf{1}.& \\
\end{eqnarray*}

\section{Concluding Remarks} 
\label{sec:concluding_remarks}
In this paper, we advance the conventional phase balancing problem to a robust look-ahead optimization framework that pursuits balanced phases in the presence of uncertainties. It is shown that imbalances among phases could be significantly reduced at the cost of a limit number of phase swapping operations. 
Many interesting directions are open for future research. For example, choosing different uncertainty sets for r-LAPB could take advantage of strong correlation among some loads. Future works also include designing approximation algorithms with optimality guarantees and exploring the benefits of controlling distributed generations \cite{xu2017adaptive,xu_adaptive_2018,zhang2018data}, electric vehicles \cite{zhang2017evaluation,zhang2018joint,chen2018impacts},  energy storage \cite{bennett_development_2015,sun_phase_2016,xu2018optimal,xu2018factoring} and demand response \cite{haider_control_2015,halder_architecture_2016,xia_energycoupon:_2017,hao_ming_scenario-based_2017}.

\section*{Acknowledgement}
This work was supported in part by Power Systems Engineering Research Center (PSERC) and in part by the Simons Institute for the Theory of Computing.
This work was partially done while the authors were visiting the Simons Institute for the Theory of Computing, UC Berkeley, and we would like to acknowledge the Gordon and Betty Moore Foundation for their generous support. 

\appendix
\section{Appendix: Equivalent Formulation of Robust Look-ahead Phase Balancing} 
Formulation (\ref{opt:robust_dynamic_phase_balancing}) is equivalent with the following:

\begin{small}
\begin{subequations}
\label{opt:robust_dynamic_phase_balancing_eq}
\allowdisplaybreaks
\begin{align}
\minimize \quad & u +\lambda v \\
\text{subject to} \quad &  h_t^\intercal p_a[t] \le u, \quad H_t^\intercal p_a[t] = a[t] - \frac{\mathbf{1}}{3}, \\
&  h_t^\intercal q_a[t] \le u, \quad H_t^\intercal q_a[t] = \frac{\mathbf{1}}{3} -a[t], \\
& h_t^\intercal p_b[t] \le u, \quad H_t^\intercal p_b[t] = b[t] - \frac{\mathbf{1}}{3}, \\
& h_t^\intercal q_b[t] \le u, \quad H_t^\intercal q_b[t] = \frac{\mathbf{1}}{3} -b[t], \\
& h_t^\intercal p_c[t] \le u, \quad H_t^\intercal p_c[t] = c[t] - \frac{\mathbf{1}}{3}, \\
& h_t^\intercal q_c[t] \le u, \quad H_t^\intercal q_c[t] = \frac{\mathbf{1}}{3} -c[t], \\
& \qquad \qquad \qquad \qquad \quad t = 1,2,\cdots,T_1, \nonumber \\
& h_t^\intercal p_a[t] \le v, \quad H_t^\intercal p_a[t] = a[T_1] - \frac{\mathbf{1}}{3}, \\
& h_t^\intercal q_a[t] \le v, \quad H_t^\intercal q_a[t] = \frac{\mathbf{1}}{3} -a[T_1], \\
& h_t^\intercal p_b[t] \le v, \quad H_t^\intercal p_b[t] = b[T_1] - \frac{\mathbf{1}}{3}, \\
& h_t^\intercal q_b[t] \le v, \quad H_t^\intercal q_b[t] = \frac{\mathbf{1}}{3} -b[T_1], \\
& h_t^\intercal p_c[t] \le v, \quad H_t^\intercal p_c[t] = c[T_1] - \frac{\mathbf{1}}{3}, \\
& h_t^\intercal q_c[t] \le v, \quad H_t^\intercal q_c[t] = \frac{\mathbf{1}}{3} -c[T_1], \\
& \qquad \qquad \quad  t = T_1+1,T_1+2,\cdots,T_2, \nonumber \\
& -w_a[t] \le a[t] - a[t-1] \le w_a[t], \\
& -w_b[t] \le b[t] - b[t-1] \le w_b[t], \\
& -w_c[t] \le c[t] - c[t-1] \le w_c[t], \\
& \sum_{t=1}^{T_1} \big( \mathbf{1}^\intercal w_a[t] + \mathbf{1}^\intercal w_b[t] + \mathbf{1}^\intercal w_c[t] \big) \le 2s, \\
& a[t] + b[t] + c[t] = \mathbf{1}, \\
& a[t], b[t], c[t] \in \{0,1\}^n, u, v \in \mathbb{R}_+,  \\
& w_a[t], w_b[t], w_c[t] \in \mathbb{R}_+^n \\
& \qquad \qquad \qquad \qquad \qquad t = 1,2,\cdots,T_1, \nonumber \\
& p_a[t], p_b[t], p_c[t], q_a[t], q_b[t], q_c[t] \in \mathbb{R}_+^k,  \\
& \qquad \qquad \qquad \qquad \qquad t = 1,2,\cdots,T_2. \nonumber
\end{align}
\end{subequations}  
\end{small}



\bibliographystyle{ieeetr}
\bibliography{myreferences}

\end{document}

%% file: hicss51-packages.tex
\usepackage[letterpaper]{geometry}
\usepackage{hicss51}
\usepackage{times}
\usepackage[none]{hyphenat}
\usepackage{url}
\usepackage{latexsym}
\usepackage{indentfirst}
\usepackage{graphicx}
\graphicspath{{images/}}
\usepackage{subfig}
\usepackage{cite}

\usepackage{amsmath,amsthm,amsfonts, amssymb}
\theoremstyle{plain}

\theoremstyle{definition}
